\documentclass[11pt,a4paper]{article}
\usepackage{amssymb, amsfonts, amsmath, latexsym}
\usepackage[english]{babel}

\begin{document}

\Large\centerline{\bf  Coarse structures on groups defined by
$T$-sequences }\vspace{6 mm}

\normalsize\centerline{\bf D. Dikranjan,  I. Protasov}\vspace{7 mm}

{\bf Abstract. } A sequence $(a_{n}) $ in an Abelian group  is called a  $T$-sequence
  if there exists a Hausdorff group topology on $G$  in which
  $(a_{n}) $ converges to $0$. For a $T$-sequence  $(a_{n}) $,  $\tau_{(a_{n}) } $
 denotes the  strongest group topology on $G$  in which
 $(a_{n}) $ converges to $0$.
The ideal  $\mathcal{I}_{(a_{n})} $ of all precompact subsets of $(G, \tau_{(a_{n}) } )$
 defines  a coarse  structure  on $G$  with base of entourages
 $\{(x, y): x-y \in P  \}$,  $P\in\mathcal{I}_{(a_{n})}. $
We prove that  $asdim \  \   (G, \mathcal{I}_{(a_{n}) }) =\infty  $ for every non-trivial
$T$-sequence  $(a_{n})$ on $G$,   and the coarse  group
  $(G, \mathcal{I}_{(a_{n}) })$ has 1 end provided that $(a_{n}) $ generates $G$.
The keypart play  asymorphic  copies of the Hamming space in
$(G, \mathcal{I}_{(a_{n})})$.

\vskip 12pt

\vskip 7pt


\vskip 15pt

Let $X$  be a set. A family $\mathcal{E}$ of subsets of $X\times X$ is called a {\it coarse structure } if
\vskip 7pt

\begin{itemize}
\item{}   each $E\in \mathcal{E}$  contains the diagonal  $\bigtriangleup _{X}$,
$\bigtriangleup _{X}= \{(x,x): x\in X\}$;
\vskip 5pt

\item{}  if  $E$, $E^{\prime} \in \mathcal{E}$ then $E\circ E^{\prime}\in\mathcal{E}$ and
$E^{-1}\in \mathcal{E}$,   where    
$$
E\circ E^{\prime}=\{(x,y): \! (\exists z \in X)((x,z) \in  E,  \  (z, y)\in E^{\prime})\} \ \mbox{ and } \ E^{-1}=\{(y,x): (x,y)\in E\};
$$
\vskip 5pt

\item{} if $E\in\mathcal{E}$ and $\bigtriangleup_{X}\subseteq E^{\prime}\subseteq E  $   then
$E^{\prime}\in \mathcal{E}$;
\vskip 5pt

\item{}  for any   $x,y\in X$, there exists $E\in \mathcal{E}$   such that $(x,y)\in E$.

\end{itemize}
\vskip 7pt

A subset $\mathcal{E}^{\prime} \subseteq \mathcal{E}$  is called a
{\it base} for $\mathcal{E}$  if, for every $E\in \mathcal{E}$, there exists
  $E^{\prime}\in \mathcal{E}^{\prime}$  such  that
  $E\subseteq E ^{\prime}$.
For $x\in X$,  $A\subseteq  X$  and
$E\in \mathcal{E}$, we write
$$
E[x]= \{y\in X: (x,y) \in E\}, \ \mbox{ and } \ E [A] = \cup_{a\in A}   \   \   E[a],
$$
 and we say that  $E[x]$ and $E[A]$
 are {\it balls of radius $E$ around} $x$  and $A$.

The pair $(X,\mathcal{E})$ is called a {\it coarse space} \cite{b9}. Alternatively, a coarse  space can be defined in terms of balls \cite{b5}, \cite{b7}.
In this case, it is called  a {\em ballean}. As pointed out in \cite{b1}, balleans and coarse structures 
as two faces of the same coin. 

\vspace{6 mm}

   A subset $A$  of a coarse space $(X, \mathcal{E})$  is called
\vspace{4 mm}

\begin{itemize}
\item{}  {\it bounded} if  $A\subseteq  E [x]$ for some $x\in X$ and $E\in \mathcal{E}$;

\item{}  {\it large} if $E [A]  = X$   for some  $E \in \mathcal{E}$.
\end{itemize}
\vskip 6pt

Given two coarse spaces $(X, \mathcal{E})$ , $(X^{\prime}, \mathcal{E}^{\prime})$
 a mapping $f: X \longrightarrow X^{\prime}$  is called
 { \it  macro-uniform}
 if, for any $E \in \mathcal{E}$, there exists
   $E^{\prime} \in \mathcal{E}^{\prime}$
    such that,  for each $x\in X$,   $f(E[x] ) \subseteq   E^{\prime}[f(x)]$.
If $f$  is a bijection such that $f$  and $f^{\prime}$  are macro-uniform   then $f$   is called an {\em asymorphism}.
 The spaces $(X, \mathcal{E})$ , $(X^{\prime}, \mathcal{E}^{\prime})$  are {\it  coarsely equivalent}  if
  $ X , X^{\prime}$ have large asymorphic  subspaces.

Now let $G$  be a group. A family $\mathcal{I}$  of subsets of $G$  is called a
 {\it group ideal}
  \cite{b6}
   if $\mathcal{I}$  contains all  finite subset and, for any
   $A, B \in \mathcal{I} $ and $C \subseteq A$,  $ \  AB^{-1} \in \mathcal{I}$
    and $C\in \mathcal{I}$.
Every group ideal  $\mathcal{I}$  defines  a coarse structure on $G$  with the base of entourages
 $\{ (x, y):  x \in  Fy\}$,  $ \ F\in \mathcal{I}$.
The  group  $G$  endowed with this coarse structure  is called a {\it coarse group},
 it is denoted by $(G, \mathcal{I})$.

We  recall that a sequence $(a_{n}) _{n \in\omega} $
 in an Abelian group  $G$ is a {\it $T$-sequence} if there exists a Hausdorff group  topology on
 $G$  in which $(a_{n}) _{n \in\omega} $  converges to $0$.
For  a  $T$-sequences
  $(a_{n}) _{n \in\omega} $ on
   $G$, we denote by  $\tau _{(a_{n})}$  the strongest group topology on $G$
in which $(a_{n}) _{n \in\omega} $
converges to $0$.
We put $A=\{ 0,  a_{n}, -a_{n}: n \in \omega\}$  and denote by
$A_{n}$  the sum on $n$  copies of $A$.

By [8, Theorem 2.3.11], $(G, \tau _{(a_{n})})$   is complete.
Hence, a subset $S$  of $G$  is totally bounded in $(G, \tau _{(a_{n})})$ if  and only if $S$  is precompact.

The ideal $\mathcal{I}_{(a_{n})}$
of all precompact subsets of $(G, \tau _{(a_{n})})$
defines a coarse structure on $G$  with the base of entourages
$\{(x,y):  x-y \in P\}$,   $P\in\mathcal{I}_{(a_{n})}$.
We denote the obtained coarse group by $(G, \mathcal{I} _{(a_{n})})$.

\vspace{6 mm}

{\bf Theorem 1}. {\it For any $T$-sequences    $(a_{n}) _{n \in\omega} $ on $G$,  the family
    $\{F+ A_{n}:  F\in [G]^{<\omega}\}$, $n\in\omega\}$
     is  a base for the ideal $\mathcal{I}_{(a_{n})}$.
         If  $G $  is generated by the set
      $\{a_{n}: n\in\omega \}$ then $\{A_{n} : n\in\omega\}$
        is a base for $\mathcal{I}_{(a_{n})}$
\vspace{5 mm}

Proof.}
Apply Lemma 2.3.2 from \cite{b8}.  \hfill  $\Box$

 \vspace{5 mm}

Given an arbitrary subset  $S$ of $G$,  the Cayley  graph Cay $(G, S)$ is a graph with the set of vertices  $G$  and the set of edges  $\{(x, y): x-y\in S\cup (-S)\}$.  \vspace{10 mm}

{\bf Theorem 2}. {\it If a
$T$-sequences    $(a_{n}) _{n \in\omega} $
 generates $G$  then the coarse group
 $(G, \mathcal{I}_{(a_{n})})$
 is asymorphic  to
Cay $(G,\{a_{n}: n\in\omega \})$.

\vspace{5 mm}

Proof.}
Apply and Theorem  1 and  Theorem   5.1.1 from \cite{b7}. \hfill  $\Box$

\vspace{5 mm}

{\bf Example. }
Let $G$  be the direct sum of groups
$\{<a_{n}>: n\in\omega\}$
 of order 2.
Clearly,
$(<a_{n}>) _{n\in\omega} $
 is a
 $T$-sequence
 on $G$.
By Theorem 1,
 the canonical bijection between
 $( G, \mathcal{I}_{(a_{n})})$
 and the Hamming space $\mathbb{H}$  of all finite subsets of $\omega$  is an asymorphism.

A  $T$-sequences    $(a_{n}) _{n \in\omega} $
is  called {\it trivial} if  $a_{n}= 0$  for all but finitely many $n\in \omega$.

\vspace{7 mm}

{\bf Theorem 3}. {\it
For any non-trivial $T$-sequences    $(a_{n}) _{n \in\omega} $
  on $G$, the coarse group  $( G, \mathcal{I}_{(a_{n})})$
  contains a subspace asymorphic to the Hamming space $\mathbb{H}$ so
  $asdim  ( G, \mathcal{I}_{(a_{n})}) =\infty$.

\vspace{7 mm}

Proof.}
Without loss of generality,  we suppose that
$\{a_{n}: n\in\omega\}$
 generates $G$  and $a_{n}\neq 0$  for each $n\in \omega$.

Given an arbitrary $T$-sequence    $(b_{n}) _{n \in\omega} $ in  $G$, we denote
   $$
   FS (b_{n}) _{n \in\omega} = \left\{\sum_{i\in F} b _{i}:  F\in [\omega]^{<\omega} \right\}
   $$
   and say that   $(b_{n}) _{n \in\omega} $
    is  FS-{\it strict} if, for any
    $H, F\in [\omega]^{<\omega}$, 
    $$
    \sum_{i\in H}  b_{i}= \sum_{i\in F}  b_{i} \  \  \  \Longrightarrow \  \  \  F=H.
    $$
We note that
$(b_{n}) _{n \in\omega}$
is  FS-strict  if, for each  $n\in\omega $,
\vspace{1mm}
$$ 
 b_{n+1}\notin \left\{\sum_{i\in F}  b_{i} -  \sum_{i\in H}  b_{i}: H, F\subseteq\{0, \ldots , n\} \right\}. \eqno(1)
$$
\vspace{1 mm}

We   assume that $(b_{n}) _{n \in\omega}$ is FS-strict and
\vspace{1 mm}

$$ 
 \mbox{if } \ b=\sum_{i\in F}  b_{i}, \    a\in A_{n} 
 \mbox{ and } b+a=\sum_{i\in H}  b_{i}   \mbox{ then } |F \triangle H|\leq n. \eqno(2)
$$
\vspace{1mm}

Then the canonical bijection $f: \mathbb{H}\longrightarrow FS (b_{n}), \  \  f(H)= \sum_{i\in H}  b_{i}$ is an asymorphism.

To construct the desired sequence  $(b_{n})_{n\in\omega}$ we rewrite  $(2)$
 in the following equivalent form
\vspace{1 mm}

 $$  
  \mbox{ if } \ i_{0} <  i_{1} <  \ldots < i_{n} <  \omega \  \mbox{ and } t _{i_{0}} ,    \ldots , t _{i_{n}} \in  \{1, -1\}
  \mbox{ then } t_{i_{0}} \  b _{i_{0}}+     \ldots +  t _{i_{n}} b _{i_{n}}\notin  A_{n}.   \eqno(3)
$$

\vspace{1mm}

We put $b_{0}= a_{0}$ and assume that  $b_{0}, \ldots , b_{n}$  have been chosen. We show how to choose  $b _{n+1}$ to satisfy $(1)$  and
\vspace{1 mm}

  $$
  \sum_{s=0}^kt _{i_{s}} b _{i_{s}}
 + t _{n+1} b _{n+1} \notin  A_{k+1}\mbox{ for }i_{0}< \ldots < i_{k}\leq n
 \mbox{ and } t_i\in\{1, -1\}\mbox{ for } i \in \{i_0,\ldots,i_k, {n+1}\}.\eqno(4) 
$$ 

  \vspace{2 mm}

We  assume that there exists a subsequence $(c_{m})_{m\in\omega}$
of $(a_{n})_{n\in\omega}$
such that $t _{i_{0}} b _{i_{0}} + \ldots +  t _{i_{k}} b _{i_{k}} +
 t  c _{m} \in  A_{k+1}$ for $t \in  \{1, -1\} $ and for each $m\in\omega$.
Every infinite subset of $A _{k+1}$  has a limit point in $A _{k}$.
Hence, $t _{i_{0}} b _{i_{0}} + \ldots +  t _{i_{k}} b _{i_{k}}  \in  A_{k+1}$
 contradicting the  choice of $b_{0}, \ldots , b_{n}$. Thus, $b _{n+1}$  can be taken from
$\{a_{m_{0}}, a_{m_{0}+1},  \ldots\}$ for some $m_{0} \in \omega$.  \hfill  $\Box$

\vspace{6 mm}

Let $(X,  \mathcal{E})$  be a coarse space. A function $f: X \longrightarrow  \{0,1\}$  is called {\it slowly oscillating} if,  for every $E\in  \mathcal{E}$ ,  there exists a bounded subset $B$ of $X$  such that $f|_{E[x]}  =  const $  for each  $x\in X\setminus   B$.
We endow $X$  with the discrete  topology,  identify the Stone-$\check{C}$ech  compactification $\beta X$  of $X$  with the set of ultrafilters on $X$   and denote
 $X^{\sharp} = \{ p\in  \beta X : $  each $P\in p$  is unbounded $\}$.
We define an equivalence $\sim$  on  $X^{\sharp}$  by the rule:  $p\sim q$
 if and only if   $f^{\beta} (p) = f^{\beta} (q) $ for every slowly oscillating
 function $f: X\longrightarrow\{0,1\} $.
The quotient $X^{\sharp} / \sim$  is called a space of  {\it ends }  or {\it binary corona}  of   $(X, \mathcal{E})$,   see  \cite{b3}, \cite{b4}.

\vspace{5 mm}

{\bf Theorem 4}. {\it If a non-trivial $T$-sequences    $(a_{n}) _{n \in\omega} $
generates $G$ then the space of  ends $( G, (a_{n})_{n\in\omega})$ is a singleton.

\vspace{5 mm}

Proof.} First we show that for every slowly oscillating function $f: G\longrightarrow\{0,1\}$ there exists an $m$ such that 
$$
f|_{G\setminus  A _{m}}  =  const .\eqno(5). 
$$ 

Indeed, by the definition of slow oscillation and Theorem 1,  there exists $m\in \omega$  such that $f|_{x+A}  =  const $ for each $x\in G\setminus A_{m}$.
We show nw that (1) holds true for this $m$.  

We take arbitrary  $y, z \in  G\setminus  A _{m} $. Since $A$ generates $G$ and contains 0, there exists 
an index $k$ such that $y,z \in A_k$, i.e., 
$$
y= b_{1} + \ldots  + b_{k} \ \ \mbox{ and }\  \ z= c_{1} + \ldots + c_{k},
$$
for appropriate 
   $ \  \  b_{1}, \ldots ,b_{k}$,   $ \ \  c_{1}, \ldots ,  c_{k} \in  A $. 
 By a property of $T$-sequences established at the end of the proof of Theorem 3,   there exists a member $a_{m_1}$ of $(a_n)$  such that
 $$ 
 a_{m_1} +  y  \notin  A _{m+1} \ \ \mbox{ and } \ \ \  a_{m_1} +  z  \notin  A _{m+1},
 $$
 since $(a_n)$ is a $T$-sequence. 
Then
$a_{m_1} +    b_{2} + \ldots  +  b_{k} \notin  A _{m}$,
   $ \ \  a_{m_1} +  c_{2}+  \ldots  +   c_{k} \notin  A _{m}$. 
   Therefore, 
$$
 f( a_{m_1} +  b_{2} + \ldots  +  b_{k}) =  f (y) \ \mbox{ and }\  \  f( a_{m_1} +  c_{2}+  \ldots  +   c_{k} )=  f(z). \eqno(6)
$$
Repeating this trick $k$   times, we can replace
$b_{2}, \ldots  ,  b_k$ and $c_{2} , \ldots  ,  c_k$, by appropriate  members $a_{m_2}, \ldots , a_{m_k}$ of $(a_n)$, as before. 
 Hence, we can replace $b_{2} + \ldots  +  b_k$ and $c_{2} + \ldots  +  c_k$, by $a_{m_2} + \ldots  +  a_{m_k}$ in (6). 
This obviously gives $f(y) = f(z)$ and proves (1). 

 Finally, the prove the assertion of the theorem, pick $p,q\in  X^{\sharp}$. In order to check that $p\sim q$ fix an arbitrary 
slowly oscillating function $f: X\longrightarrow\{0,1\} $. We have prove that $f^{\beta} (p) = f^{\beta} (q) $. Pick an $m$ with (1). Since $G\setminus  A _{m}\in p \cap q$, for every $P\in p$ and for every $Q\in q$ we have     $P_1:= P \setminus A _{m}\in p $ and $Q_1:= Q\setminus  A _{m}\in q$
 and $f|_{P_1} = f|_{Q_1}$ is constant in view of (5). This proves that     $f^{\beta} (p) = f^{\beta} (q) $. 
 $ \Box$
\vspace{6 mm}

{\bf Remark 1.}
An $FS$-strict sequence generating $G$  is called  a  {\it base} for $G$.
By [2, Theorem 3.3], a countable group $G$  has a base if and only if $G$  has no elements of odd order.
If $(a_{n}) _{n \in\omega} $ is a base  of $G$  then we have the natural  bijection
 $f:  \mathbb{H} \longrightarrow G$,  $  \  \  f(F) = \sum _{i\in F}  a_{i}$.
If  in addition $(a_{n}) _{n \in\omega} $
is a  $T$-sequence,
 we cannot state that  $f$ is an asymorphism  between $\mathbb{H}$  and $(G, \mathcal{I}  _{(a_{n}) } )$.
For example, we may take the base
$((-1)^{n}  2^{n} )  _{n\in \omega}$ of $\mathbb{Z}$.
\vspace{6 mm}

We denote by $\mathbb{H}(n)$ the direct sum of $\omega$  copies of
 the cyclic group  of  order  $n$  and endow
  $\mathbb{H}(n)$  with the Hamming metric  $h$, where for $a, b\in \mathbb{H}(n)$
  $h (a, b)$ is the  number  of distinct  coordinates of $a$  and $b$. If $n,  m$  have the same set of prime  divisors  that
 $\mathbb{H}(n)$  and $\mathbb{H}(m)$ are  asymorphic.
\vspace{5 mm}

{\bf Question 1.} {\it
When precisely $\mathbb{H}(n)$  and $\mathbb{H}(m)$   are asymorphic?  Coarsely equivalent?}
\vspace{6 mm}


{\bf Remark 2.}
Given a group ideal  $\mathcal{I}$ on $G$,  the coarse
 group $(G, \mathcal{I})$  has  bounded geometry if exists $F \in  \mathcal{I}$  such that,  for each  $H\in \mathcal{I}$,  we have $H\subseteq FK$  for some  $K\in [G] ^{<\omega}$.
For each  non-trivial
$T$-sequence
 $(a_{n})_{n\in\omega}$ on $G$,
 the coarse group $(G, \mathcal{I} _{(a_{n}))}$  is not of bounded geometry.
\vspace{5 mm}

{\bf Remark 3.}
We say  that a subset $S$  of a topological group  $(G, \tau)$ is  $R${\it -bounded} (see \cite{b10}) if, for any  neighbourhood  $U$ of the identity of $G$,  there exist
$K\in [G]^{< \omega}$  and $n\in \omega$ such that
$S\subseteq (UK)^{n}$.
The  family $\mathfrak{R}_{\tau}$  of all  $R$-bounded  subsets of $(G, \tau)$  is a group  ideal.
For any  $T$-sequence on  $G$,  we have $\mathfrak{R}_{\tau(a_{n})}= \mathcal{I}_{a_{n}}$ .
\vspace{5 mm}

{\bf Question 2.} {\it
Let $G$  be a non-discrete countable  metrizable Abelian group and let
 $\mathcal{I} _{seq}$   denotes the smallest group ideal  on $G$  containing all converging sequences.
 Is $asdim \  (G,  \mathcal{I}_{seq}) = \infty$ ?}

\end{document}